%%%%%%%%%%%%%%%%%%%%%%%%%%%%%%%%%%%%%%%%%%%%%%%

%%%%%%%%%%%%%%%%%%%%%%%%%%%%%%%%%

\documentstyle{amsppt}\magnification=1200\loadbold

\vsize 22 true cm \hsize 16 true cm

\parskip=5.0pt

\topmatter

\title Regularity and relaxed problems of minimizing biharmonic maps into spheres
\endtitle

\rightheadtext{biharmonic maps} \leftheadtext{M.-C. Hong and C.Y.
Wang}

\author{Min-Chun Hong  and Changyou Wang }\endauthor

\address {Min-Chun Hong, Department of Mathematics, University of
Queensland, Brisbane,   QLD 4072 , Australia
\medskip   Changyou Wang, Department of Mathematics, University of
Kentucky, Lexington, KY 40506, USA}
\endaddress

%\subjclass{AMS  }\endsubjclass

\abstract {For $n\ge 5$ and $k\ge 4$, we show that any minimizing
biharmonic map from $\Omega\subset \Bbb R^n$ to $S^k$ is smooth
off a closed set whose Hausdorff dimension is at most $n-5$. When
$n=5$ and $k=4$, for a parameter $\lambda\in [0,1]$ we introduce a
$\lambda$-relaxed energy $\Bbb H_{\lambda}$ of the Hessian energy
for maps in $W^{2,2}(\Omega ; S^4)$ so that each minimizer
$u_{\lambda}$ of $\Bbb H_{\lambda}$ is also a biharmonic map. We
also establish the existence and partial regularity of a minimizer
of $\Bbb H_{\lambda}$ for $\lambda \in [0,1)$. }
\endabstract
\endtopmatter

\document

\def\R{{\Bbb R}}

  \def\laplacian{\bigtriangleup}

  \def\a{\alpha}

\def\intave#1{-\kern-10.7pt\int_{\,#1}}

\def\b{\beta}

\def\D{\nabla}

\def\<{\langle}

\def\>{\rangle}

\def\({\left(}

\def\){\right)}

\def\limsup{\operatornamewithlimits{lim\,sup}}

\def\intave#1{-\kern-10.7pt\int_{\,#1}}

\head {\bf 1. Introduction}\endhead

For $n\ge 5$ and $k\ge 4$, let $\Omega\subset\R^n$ be a bounded
smooth domain and $S^k\subset \R^{k+1}$ be the unit sphere. Define
$$W^{2,2}(\Omega, S^k)=\{u\in W^{2,2}(\Omega, \R^{k+1})|
\quad  |u(x)|=1 \text { for a.e. } x\in \Omega\}$$
and also
define, for a given map $u_0\in W^{2,2}(\Omega, S^k)$,
$$W_{u_0}^{2,2}(\Omega, S^k)=\{u\in W^{2,2}(\Omega, S^k)
|\quad u-u_0|_{\partial\Omega}= \nabla (u-u_0)|_{\partial\Omega}=0
\hbox{ in the trace sense }\}.$$

The hessian energy functional on $W^{2,2}(\Omega, S^k)$ is defined
by
$$\Bbb H(u)=\int_{\Omega} |\laplacian u|^2\,dx, \ \
\forall u\in W^{2,2}(\Omega,S^k).\tag 1.1$$

Recall that a map $u\in W^{2,2}(\Omega,S^k)$ is a (weakly)
biharmonic map  if $u$ is a critical point of $\Bbb H(\cdot )$ in
$W^{2,2}(\Omega, S^k)$ so that it satisfies the Euler-Lagrange
equation:
$$-\laplacian^2 u=( |\laplacian u|^2 +2\nabla\cdot(\nabla u\cdot\laplacian u)
-\laplacian |\nabla u|^2)u\tag 1.2$$ in the distribution sense,
where $\nabla\cdot$ is the divergence operator in $\R^n$ and
$\cdot$ is the inner product in $\R^{k+1}$.

A typical class of biharmonic maps is given by minimizing
biharmonic maps. A map $u\in W^{2,2}(\Omega, S^k)$ is a minimizing
biharmonic map  if it satisfies
$$\Bbb H (u)\leq \Bbb H(w), \ \ \forall w\in W^{2,2}_{u}(\Omega, S^k).\tag 1.3$$

The study of minimizing biharmonic maps into spheres was initiated
by Hardt-Mou \cite{HM}. Chang-Wang-Yang \cite {CWY} established
the partial regularity for weakly stationary biharmonic maps into
spheres, ie. if $u\in W^{2,2}(\Omega ;S^k)$ be a weakly stationary
biharmonic map, then $u\in C^\infty(\Omega,S^k)$ for $n=4$  and
$u\in C^\infty(\Omega\setminus\Sigma,S^k)$ with $\Cal
H^{n-4}(\Sigma)=0$ for $n\ge 5$. Very recently, the main theorems
of \cite{CWY} have been generalized by the second author
\cite{W1,2,3} for stationary biharmonic maps into any compact
smooth Riemannian submanifold $N$ of  the Euclidean spaces. The
results in \cite {CWY}, \cite{W1,2,3} give the partial regularity
theorem for minimizing biharmonic maps for $n\ge 5$ since a
minimizing biharmonic map  is stationary. However the minimality
of the biharmonic map $\Phi(x,y)={x\over |x|}: B^5\times B^{n-5}
\to S^4$ (see Proposition A1 of \S5 below) indicates that the
dimension of the singular set $\Sigma$ for minimizing biharmonic
maps may be smaller than $n-4$. Compared with the optimal partial
regularity for minimzing harmonic maps by Giaquinta-Giusti \cite
{GG} and by Schoen-Uhlenbeck \cite {SU}, it is natural to ask the
following question:

\noindent{\it Does the singular set of a minimizing biharmonic map
have Hausdorff dimension at most $n-5$?}

In this aspect, we have
\proclaim {Theorem A} For $n\ge 5$ and
$k\ge 4$, let $u\in W^{2,2}(\Omega, S^k)$ be a minimizing
biharmonic map and denote by ${\Cal S}(u)$ the singular set of
$u$. Then ${\Cal S}(u)$ is discrete for $n=5$, and has its
Hausdorff dimension at most $n-5$ for $n\ge 6$.
\endproclaim

One of key ingredients to prove Theorem A is  to derive  an extension
inequality (2.1) for the hessian energy  for maps into $S^k$ with
$k\ge 4$, inspired by Hardt-Lin's extension Lemma
(see \cite {HL}, \cite {HKL}) in the context of harmonic maps into
a simply connected manifold. An important consequence of (2.1) is
to obtain the Caccioppoli inequality (2.7) for miminimizing
biharmonic maps so that {\it a weakly convergent sequence of
minimizing biharmonic maps in $W^{2,2}(\Omega,S^k)$ converges
strongly in $W^{2,2}_{\hbox{loc}}$ to a minimizing biharmonic map.}
Combined this fact with the energy monotonicity inequality
for stationary biharmonic maps due to \cite{CWY}, it guarantees
that a refinement of Federer dimension reduction scheme \cite{F}
is applicable so that Theorem A follows from a similar argument as
one in \cite {S}.

For $n=5$, it follows from Theorem A that suitable rescalings at
each singular point of a minimizing biharmonic map $u\in
W^{2,2}(\Omega,S^4)$ yields a minimizing biharmonic map of the
form $\Psi({x\over |x|})$ for some $\Psi\in C^\infty(S^4,S^4)$.
Motivated by the work by Brezis-Coron-Lieb \cite{BCL} on
minimizing harmonic maps from $B^3$ into $S^2$, it will be an
interesting problem to study the map $\Psi$ given as above.
Inspired by the problem proposed by Hardt-Lin \cite{HL1} on the
context of harmonic maps from $B^3$ to $S^2$, we would also like
to ask the following question:

\noindent {\it For any given map $\psi\in C^\infty(S^4,S^4)$ with
zero degree,  is there  a biharmonic map $u\in
C^\infty(\overline{B^5},S^4)$ with $u=\psi$ on $\partial B^5$?}

The example A2 of Section 5 indicates that any minimizing
biharmonic map extension of some boundary map with degree zero has
singularities. In order to study this problem, we extend the idea
of a relaxation of the Dirichlet energy functional of harmonic
maps from $B^3$ to $S^2$ by Bethuel-Brezis-Coron \cite {BBC} and
Giaquinta-Modica-Souc\`ek \cite{GMS}. More precisely, we hope to
introduce a relaxed energy functional for biharmonic maps from
$\Omega\subset \R^5$ to $S^4$.

For a $W^{2.2}$-map from $\Omega\subset \R^5$ to $S^4$, the
D-field of $u$, $D(u)=(D_1(u),...,D_5(u))\in L^1(\Omega, \R^5)$,
is defined by
$$D_1(u)=\text {det }\left (u,\frac {\partial u}{\partial x_2}, ...,  \frac {\partial
u}{\partial x_5}\right ),\quad...,\quad  \,D_5(u)=\text {det
}\left (u, \frac {\partial u}{\partial x_1}, \cdots,  \frac
{\partial u}{\partial x_4}\right ).$$

For any given $\phi\in W^{2,2}(\Omega,S^4)\cap
C^\infty(\partial\Omega,S^4)$, define
$$L(u):=\frac 1 {\sigma_4}\sup_{\xi :\Omega\to \R, \|\D\xi\|_{L^{\infty}}\leq 1}
\{\int_{\Omega} D(u)\cdot \D\xi\,dx-\int_{\partial\Omega}D(u)\cdot
\nu\xi\,dH^{n-1}\}, \ \forall u\in W^{2,2}_\phi(\Omega,S^4).
$$

For any $\lambda\in (0,1]$, the  $\lambda$-relaxed hessian energy
functional is  defined by
$$\Bbb H_{\lambda}(u):=\Bbb H(u) +16 \lambda \sigma_4 L(u),
\ \forall u \in W^{2,2}_\phi(\Omega,S^4).\tag 1.3$$

Throughout the paper, we denote by $\sigma_k=\Cal H^k(S^k)$ the
area of the unit sphere $S^k\subset \R^{k+1}$ for $k\ge 4$.

Then we have

\proclaim {Theorem B}

(a) For any $\lambda\in (0,1]$, $\Bbb H_{\lambda}$ is sequentially
lower semi-continuous in $W^{2,2}_\phi(\Omega,S^4)$ for the weak
$W^{2,2}$-topology,

(b) for any $\lambda \in (0,1]$, there exists a $u_{\lambda}\in
W^{2,2}_\phi(\Omega,S^4)$ which minimizes $\Bbb
H_{\lambda}(\cdot)$ among $W^{2,2}_{\phi}(\Omega , S^4)$-maps, and

(c)for any $\lambda\in (0,1)$, $u_\lambda$ is a weakly biharmonic
map satisfying $u_{\lambda}\in
C^\infty(\Omega\setminus\Sigma_\lambda,S^4)$, with $\Cal
H^{1-\delta}(\Sigma_\lambda)=0$ for some $\delta>0$.
\endproclaim

Finally,   modifying the arguments of \cite {BCL}, there exist
infinitely many weak biharmonic maps in $W^{2,2}_{x/|x|}(\Omega ;
S^4)$. It will be an interesting question to establish this result
for general boundary data. To do it, one needs to establish the
boundary regularity of a minimizing harmonic maps, but this is
unknown. The partial regularity has been established in \cite
{LW}.

The paper is organized as follows. In Section 2, we derive a
Caccioppoli's inequality for $Q$-minimizing biharmonic maps. In
Section 3, we prove a partial regularity for $Q$-minimizing
biharmonic maps and also present a proof of Theorem A. In Section
4, we prove Theorems B. In Section 5 is an appendix and several
elementary facts will be given.

\subhead Acknowledgment  \endsubhead  We would like to thank
Professor Mariano Giaquinta for useful comments and discussions.
The research of the first author is supported by the Australian
Research Council.

\head {\bf 2. Caccioppoli's inequality}\endhead

In this section, we consider, for any $Q\ge 1$, $Q$-minimizing
biharmonic maps from $\R^n$ to $S^k$ for $n\ge 5$ and $k\ge 4$,
and establish the Cacciopoli inequality. The idea is inspired by
the Hardt-Lin's extension Lemma (see \cite{HL}, \cite{HKL}).

\proclaim {Definition 2.0}($Q$-minimizing biharmonic map) Let $Q$
be a constant with $1\le Q<\infty$. For $n\ge 5$ and $k\ge 4$, a
map $u\in W^{2,2}(\Omega,S^k)$ is called a $Q$-minimizing
biharmonic map if (i) $u$ is a weakly biharmonic map  and (ii) $u$
satisfies
$$\Bbb H(u)\le Q \Bbb H(v), \ \forall v\in W^{2,2}_u(\Omega,S^k).\tag 2.0$$
\endproclaim

It is clear that any minimizing biharmonic map is a
$Q$-minimizing biharmonic map with $Q=1$. Now we have

\proclaim {Lemma 2.1}(Extension Lemma) For $n\ge 5$ and $k\ge 4$, let $\Omega\subset \R^n$
 be any bounded smooth domain. Then for any map $v\in W^{2,2}(\Omega,\R^{k+1})$
with $|v(x)|=1$ and $\nabla v(x)\in T_{v(x)}S^k$ for a.e.
$x\in\partial\Omega$, there exists a map $w\in W^{2,2}(\Omega,
S^k)$ such that $w=v$,$\nabla w=\nabla v$ on $\partial\Omega$, and

$$\int_{\Omega} |\laplacian w|^2\,dx\leq C\int_{\Omega }(| \laplacian
v|^2+|\nabla v|^4)\,dx,\tag 2.1$$
where $C>0$ is independent of
$u$ and $\Omega$.
\endproclaim

\demo {Proof} For any $a\in \R^{k+1}$ with $|a|\leq \frac 12 $,
consider the map

$$w_a(x)=\frac {v(x)-a}{|v(x)-a|},\quad x\in \Omega.$$

Then a simple calculation gives

$$\nabla w_{a}(x)=|v(x)-a|^{-1} \nabla v(x)-|v(x)-a|^{-3}(v(x)-a)\otimes
(v(x)-a)\nabla v(x), \ \forall x\in\Omega,$$

we have, by taking one more derivative of this identity,

$$ |\laplacian w_a(x)|\leq C\left (\frac {| \laplacian v(x)|}{|v(x)-a|}+ \frac {|\nabla
v(x)|^2 }{|v(x)-a|^2}\right ), \ \forall x\in\Omega. \tag 2.2$$

Integrating (2.2) over $\Omega\times B_{1\over 2}^{k+1}$ with
respect to $(x,a)$  and applying Fubini's theorem, we have

$$\split&\int_{B_{1/2}^{k+1}}\int_{\Omega} | \laplacian w_a|^2(x)\,dx\,da
=\int_{\Omega}\int_{B_{1/2}^{k+1}}| \laplacian w_a|^2(x)\,dx\,da\\
&\leq C \int_{\Omega}\(|\laplacian v|^2+|\nabla
v|^4\)(x)\cdot\left
 [\int_{B_{1\over 2}^{k+1}}
\({1\over |v(x)-a|^2}+{1\over |v(x)-a|^4}\)\,da\right ]\,dx\\
&\leq C\int_{\Omega}\(|\nabla v|^4+|\laplacian
v|^2\)(x)\,dx,\endsplit$$ where we use the fact that
$$ \int_{B_{1/2}^{k+1}} {1\over |v(x)-a|^4}\,da\leq
c(k)=\cases\frac {16} 9 \frac {\sigma_k}{k+1}, &\text {for $|v|\geq 1$}\\
\frac 1{k-3}(\frac 32 )^{k-3} \sigma_k, & \text {for $|v|\leq 1$}.
\endcases $$

Therefore we can find an $a_0\in B_{1\over 2}^{k+1}$ such that
$$\int_{\Omega} | \laplacian w_{a_0}|^2(x)\,dx\leq C\int_{\Omega} (| \laplacian
v|^2+|\nabla v|^4)(x)\,dx.\tag 2.3$$
For $a\in B_{1\over
2}^{k+1}$, define
$$\Pi_a(\xi)=\frac {\xi -a}{|\xi -a|}: S^k\to S^k.$$
It is easy to see that $\Pi_a$ is a $C^2$ diffeomorphism of $S^k$
onto itself. In fact,
$$\Pi^{-1}_a(\xi)=a+[ (a\cdot\xi )^2+(1-|a|^2)]^{1/2}\xi, \ \forall \xi\in S^k.
$$
In particular, we have
$$\max_{a\in B_{1\over 2}^{k+1}}\(\|\nabla \Pi_a^{-1}\|_{C^0(S^k)}+
\|\nabla^2 \Pi_a^{-1}\|_{C^0(S^k)}\)=\Lambda<\infty.$$ Now we set
$$w(x)=\Pi^{-1}_{a_0}\circ w_{a_0}(x)=\Pi^{-1}_{a_0}\circ\Pi_{a_0}(v(x)),
\ ,\forall x\in \Omega.$$

Note
$$\split\nabla w(x)&=\nabla\Pi^{-1}_{a_0}(w_{a_0}(x))\nabla w_{a_0}(x),\\
\nabla^2 w(x)&=\nabla\Pi_{a_0}^{-1}(w_{a_0}(x))\nabla^2 w_{a_0}(x)
+\nabla^2 \Pi_{a_0}^{-1}(w_{a_0}(x))(\nabla w_{a_0}(x), \nabla
w_{a_0}(x)) \endsplit$$ and
$$|\nabla w_{a_0}(x)|^2\le |\laplacian w_{a_0}(x)|, \text{ for a.e. }x\in\Omega$$
due to the fact that $|w_{a_0}(x)|=1$ for a.e. $x\in\Omega$. Then
we have
$$|\laplacian w|(x)\leq C(\Lambda)[|\laplacian w_{a_0}|+|\nabla w_{a_0}|^2](x)
\le C(\Lambda)|\laplacian w_{a_0}|(x), \, \forall x\in\Omega.\tag
2.4$$

Combined this with (2.3), it implies that $w\in W^{2,2}(\Omega,
S^k)$ and satisfies (2.1). To see that $w$ has the same trace as
$v$ on $\partial\Omega$, observe that, since $\Pi_{a_0}^{-1}\circ
\Pi_{a_0}|_{S^k}$ is the identity map, $w=v$ on $\partial\Omega$.
Moreover, since for any $x\in\partial\Omega$ we have $\nabla
w(x)=\nabla(\Pi_{a_0}^{-1}\circ\Pi_{a_0})(v(x))(\nabla v(x))$,
$\nabla v(x)\in T_{v(x)}S^k$, and $\nabla(\Pi_{a_0}^{-1}\circ
\Pi_{a_0})(v(x)): T_{v(x)}S^k\to T_{v(x)}S^k$ is the identity map,
we have $\nabla w=\nabla v$ on $\partial\Omega$. The proof of
Lemma 2.1 is complete.
 \qed
\enddemo

We follow the iteration method in \cite{G} to get
\proclaim {Lemma
2.2} For $0\leq r_0<r_1<\infty$, let $f: [r_0,r_1]\to (0,\infty)$
be a mesurable function. Suppose that there exist $\theta\in
(0,1)$, $A>0$ $B>0$, $\alpha$, and $\beta>0$ such that for
$r_0\leq t<s\leq r_1$ we have
$$f(t)\leq \theta f(s)+[A(s-t)^{-\a} +B(s-t)^{-\b}].\tag 2.5$$
Then for all $r_0\leq \rho R\leq r_1$ we have
$$ f(\rho )\leq C[A(R-\rho )^{-\a} + B(R-\rho )^{-\b}],\tag 2.6$$
where $C=C(\a,\b,\theta)>0$.
\endproclaim

Now we have

 \proclaim {Lemma 2.3} (Cacciopoli's
inequality) For $1\leq Q<\infty$, $n\geq 5$ and $k\geq 4$, let $u$
be a $Q$-minimizer of $\Bbb H$ in $W^{2,2}_{u_0} (\Omega , S^k)$.
Then for all $x_0\in\Omega$ and all $R<\text{dist} (x_0,
\partial\Omega )$, we have
$$\int_{B_{R/2}(x_0)}|\laplacian u|^2\,dx\leq C R^{-4}\int_{B_R(x_0)}( | u-u_{x_0,
R}|^2+| u-u_{x_0, R}|^4 )\,dx \tag 2.7$$ for some constant $C>0$,
where $u_{x_0,R}=\frac 1{|B_R(x_0)|}\int_{B_R(x_0)} u\,dx$ is the
average of $u$ over $B_R(x_0)$.
\endproclaim
\demo {Proof} It follows from Lemma 2.1 and the Q-minimality of
$u$  that
$$\int_{\Omega}|\laplacian u|^2\,dx\leq CQ\int_{\Omega}\(|\laplacian v|^2+ |\nabla v|^4\)\,dx\tag
2.8
$$
for any $v\in W^2_{u}(\Omega; \R^{k+1})$, where $C$ is a positive
constant.

For any $s$, $t$ with $R/2\leq t<s\leq \frac {3R}4$, let $\phi$ be
a cut-off function in $B_s(x_0)$ such that $0\leq \phi \leq 1$ in
$B_s(x_0)$, $\phi\equiv 1$ in $B_t(x_0)$, $\phi\equiv 0$ outside
$B_s(x_0)$, $|\nabla \phi|\leq C(s-t)^{-1}$ and $|\laplacian \phi
|\leq CR^{-2}$, where $C$ is a constant independent of $s$, $t$
and $R$. Taking $v(x)=u(x)-\phi [u(x)-u_{x_0,R}]$, we have
$$\nabla v=(1- \phi )\nabla u-\nabla\phi [u(x)-u_{x_0,R}] $$ and
$$\laplacian v=(1-\phi )\laplacian u-\laplacian \phi [u(x)-u_{x_0,R}]-2\nabla \phi \nabla u.$$

By (2.8), we obtain
$$\split
\int_{B_t}|\laplacian u|^2\,dx&\leq C\int_{B_s\backslash B_t}(
|\laplacian u|^2+|\nabla
u|^4)\,dx   + \frac {C}{(s-t)^2}\int_{B_s}|\nabla u|^2 \,dx\\
\quad &+\frac C{(s-t)^4}\int_{B_s} |u-u_{x_0,R}|^2\,dx+\frac
C{(s-t)^4}\int_{B_s} |u-u_{x_0,R}|^4\,dx\endsplit \tag 2.9$$ for
all $s$, $t$ with $\frac R2\leq t<s\leq R$.

Noticing $|u|=1$, we have
$$|\nabla u|^2\leq |\laplacian u|.\tag 2.10$$

By the filling hole trick in (2.9), there exists a positive
$\theta <1$ such that
$$\split
\int_{B_t}|\laplacian u|^2\,dx&\leq \theta \int_{B_s}
|\laplacian u|^2\,dx +\frac  C{(s-t)^2}\int_{B_s} |\nabla u|^2\,dx\\
\quad &+\frac C{(s-t)^4}\int_{B_s} |u-u_{x_0,R}|^2\,dx+\frac
C{(s-t)^4}\int_{B_s} |u-u_{x_0,R}|^4\,dx\endsplit
$$
for $\frac R2 \leq t<s\leq \frac {3R}4$. Then it follows from
Lemma 2 to obtain
$$\split
\int_{B_{R/2}(x_0)}|\laplacian u|^2\,dx\leq &C
R^{-2}\int_{B_{\frac {3R}4}(x_0)}|\nabla u|^2\,dx+R^{-4}\int_{B_{R}(x_0)}|  u-(u)_{x_0, R}|^2\,dx\\
&+R^{-4}\int_{B_R(x_0)}|  u-(u)_{x_0, R}|^4\,dx.
\endsplit \tag 2.11
$$

Let $\phi \in C^{\infty}_0(B_R(x_0))$ be a cut-off function with
$\phi \equiv 1$ in $B_{\frac {3R}4 }(x_0)$, $0\leq \phi \leq 1$
and $|\nabla \phi \leq \frac CR$. Integrating by parts, we have
$$\int_{B_R(x_0)} \phi^2 |\nabla u|^2\,dx =-\int_{B_R(x_0)}
\laplacian u \cdot (u-u_{x_0, R})\phi^2 -2\int_{B_R(x_0)} \nabla
u\cdot (u-u_{x_0, R})\phi\nabla\phi\,dx.$$ Then
$$\int_{B_{\frac {3R}4 }(x_0) } |\nabla u|^2\,dx\leq \varepsilon R^2\int_{B_R(x_0)}
|\laplacian u|^2\,dx +\frac C{R^2} \int_{B_R(x_0)} |u-u_{x_0,
R}|^2\,dx\tag 2.12$$ for a sufficiently small $\varepsilon$.

By (2.11)-(2.12) with a sufficiently small $\varepsilon$,  the
claim (2.7) follows from the standard trick (e.g. \cite {S; Lemma
2 of Chapter 2}).\qed
\enddemo

As a direct consequence of Lemma 2.3, we have the following
reverse H\"older inequality for $Q$-minimizing biharmonic maps.

\proclaim {Proposition 2.4} Let $1\leq Q<\infty$, $n\geq 5$ and
$k\geq 1$. Suppose that $u$ be a Q-minimizing biharmonic map in
$W^{2,2}(\Omega ; S^k)$. Then there exists an exponent $p>4$ such
that $u\in W^{2,p}_{loc} (\Omega , \R^5)$. Moreover, for all
$x_0\in\Omega$ and $R<\text{dist}(x_0,\partial\Omega)$, we have
$$\left (\intave {B_{R/2}(x_0)|} (|\laplacian u|^2+1)^{p/2}\,dx\right )^{1/p}\leq C \left (\intave {B_R (x_0)}
(|\laplacian u|^2+1)\,dx\right )^{1/2},\tag 2.13
$$
for some constant $C$ depending $n$, $k$ and $Q$, where we denote
by the  average integration over $B_R(x_0)$
$$\intave
{B_{R}(x_0)}f(x)\,dx=\frac 1{|B_R(x_0)}\int_{B_R(x_0)}f(x)\,dx.$$
\endproclaim
\demo {Proof} By the Poincare inequality, we have
$$\int_{B_R(x_0)}|u-u_{x_0,R}|^4\,dx\leq CR^{4+n(1-\frac 4q )} \left ( \int_{B_R(x_0) }|\nabla u|^q\right
)^{4/q}
$$
and
$$\int_{B_R(x_0)}|u-u_{x_0,R}|^2\,dx\leq CR^{2+n(1-\frac 2q )} \left ( \int_{B_R(x_0) }|\nabla u|^q\right
)^{2/q}
$$ for some $q<4$. Then it follows from  (2.7) that
$$\intave {B_{R/2}(x_0)} |\laplacian u|^2 \,dx\leq C \left (\intave {B_{R}(x_0)} |\nabla u|^q \,dx\right
)^{4/q}+ C \left (\intave {B_{R}(x_0)} |\nabla u|^q \,dx\right
)^{2/q}$$ for some $q<4$. This implies from (2.7) that
$$\intave {B_{R/2}(x_0)}\left [ |\laplacian u|^2 +1\right
]\,dx\leq C \left (\intave {B_{R}(x_0)}  \left [ |\laplacian
u|^2+1\right ]^{q/4} \,dx\right )^{4/q}$$ for every
$B_{R}(x_0)\subset\Omega$ with $q<4$. By the reverse H\"older
inequality (cf. \cite {G; page 122}), there exists an exponent
$p>2$ such that for every $B_R(x_0)\subset \Omega$
$$\left (\intave {B_{R/2}(x_0)} \left [ |\laplacian u|^2 +1\right
]^{p/2}\,dx\right )^{2/p}\leq C \intave {B_{R}(x_0)}  \left [
|\laplacian u|^2+1\right ] \,dx.$$ By the  standard $L^p$-local
estimate of linear elliptic equations of second order, $u\in
W_{loc}^{2,p}(\Omega,\R^5)$. \qed

\enddemo

\head {\bf 3. Proof of Theorem A}\endhead

In this section, we  present a proof of Theorem A. The proof
consists of three steps: (i)regularity under the smallness of
renormalized hessian energy, (ii)$W^{2,2}$-compactness of the
space of minimizing biharmonic maps  and (iii) blow-up argument
utilizing both the energy monotonicity inequality  (\cite{CWY})
and Federer's dimension  deduction  \cite{F}.

\proclaim {Lemma 3.1} For any $1\le Q<\infty$, $n\ge 5$, and $k\ge
4$, there exists an $\epsilon_0 \in (0,1)$ such that if $u\in
W^{2,2}(\Omega,S^k)$ is an $Q$-minimizing biharmonic map
satisfying
$$ (2R)^{4-n}\int_{B_{2R}(x_0)}|\laplacian u|^2\le \epsilon_0^2,
\text{ for }B_{2R}(x_0)\subset\Omega \tag 3.1$$ then $u\in
C^{\infty}(B_R(x_0),S^k)$ and satisfies
$$\|u\|_{C^{l}(B_R(x_0)}\le C(n,k,Q,\epsilon_0,l), \ \forall l\ge 1.\tag 3.2$$
\endproclaim

\demo{Proof} It is based on the following decay estimate: there
exists an $\theta_0\in (0,1)$ such that
$$(\theta_0R)^{4-n}\int_{B_{\theta_0R}(x_0)}|\laplacian u|^2\,dx
\le ({1\over 2}) R^{4-n}\int_{B_R(x_0)}|\laplacian u|^2\,dx. \tag
3.3$$

Once (3.3) is established,   the regularity of $u$ follows from
the standard iterations and suitable applications of Morrey's
Lemma \cite {M}. To prove (3.3), we argue by contradiction as
follows (see \cite{HL} for similar arguments for harmonic maps).
By rescalings, we may assume that $x_0=0$ and $R=1$. Suppose (3.3)
is false. Then there exists a sequence of minimizing biharmonic
maps $\{u_i\}\subset W^{2,2}(B_1,S^k)$ such that
$$ \int_{B_1}|\laplacian u_i|^2\,dx =\epsilon_i^2\rightarrow 0$$
but we have, for any $\theta\in (0,1)$,
$$\theta^{4-n}\int_{B_\theta}|\laplacian u_i|^2\,dx>{1\over 2}
\int_{B_1}|\laplacian u_i|^2\,dx=\frac 12 \epsilon_i^2. \tag 3.4$$
Define the blow-up sequence $v_i(x)={u_i(x)-(u_i)_{1}\over
\epsilon_n}: B_{1}\to \R^{k+1}$. It is easy see
$$(v_i)_{1}=0, \quad
\int_{B_1}|\laplacian v_i|^2\,dx=1, \quad \int_{B_1}|\nabla
v_i|^4\,dx\le 1.$$
Therefore we may assume that $v_i\rightarrow
v_\infty$ weakly in $W^{2,2}(B_{1})$, strongly in $W^{1,2}(B_1)$
and $L^4(B_1)$. Since $u_i$ satisfies (1.2), it is easy to see
that $v_i$ satisfies
$$-\laplacian^2 v_i=\epsilon_i(|\laplacian v_i|^2+2\nabla\cdot(\nabla v_i\cdot\laplacian v_i)
-\laplacian |\nabla v_i|^2)u_i, \ \hbox{ in }B_1. \tag 3.5$$
Letting $i$ tend to infinity, we see $\int_{B_1}|\laplacian
v_\infty|^2\le 1$ and
$$\laplacian^2 v_\infty =0, \ \hbox{ in }B_1.$$
Therefore $v_\infty\in C^\infty(B_1,\R^{k+1})$ and satisfies, for
any $\theta\in (0,{1\over 2})$,
$$\theta^{2-n}\int_{B_\theta}|\nabla v_\infty|^2\,dx\le C\theta^2, \ \
\theta^{-n}\int_{B_\theta}|v_\infty-(v_\infty)_\theta|^4\,dx\le
C\theta^4. \tag 3.6$$
Therefore, for $i$ sufficiently large, we
have
$$\theta^{2-n}\int_{B_\theta}|\nabla u_i|^2\,dx\le C\theta^2\epsilon_i^2, \ \
\theta^{-n}\int_{B_\theta}|u_i-(u_i)_\theta|^4\,dx\le
C\theta^4\epsilon_i^4. \tag 3.7$$ This, combined with (2.7) of
Lemma 2.3, implies that for any $\theta\in (0,{1\over 4})$
$$\theta^{4-n}\int_{B_\theta}|\laplacian u_i|^2\,dx
\le C(\theta^2+\theta^4\epsilon_i^2)\epsilon_i^2\le
C\theta^2\epsilon_i^2.$$ This contradicts with (3.4), if we choose
sufficiently small $\theta\in (0,{1\over 4})$. \qed
\enddemo

As a consequence of Lemma 3.1, we have the following partial
regularity for $Q$-minimizing biharmonic maps.

\proclaim{Corollary 3.2} For any $1\le Q<\infty$, $n\ge 5$, and
$k\ge 4$. Suppose that $u\in W^{2,2}(\Omega,S^k)$ is an
$Q$-minimizing biharmonic map. Then there exist a closed set
$\Sigma\subset\Omega$ and an $\delta>0$ such that $u\in
C^\infty(\Omega\setminus\Sigma,S^k)$ and ${\Cal
H}^{n-4-\delta}(\Sigma)=0$.
\endproclaim
\demo{Proof} It follows from Lemma 3.1 that the singular set of
$u$ is given by
$$\Sigma=\{x\in\Omega | \ \liminf_{r\rightarrow 0}
r^{4-n}\int_{B_r(x)}|\laplacian u|^2\,dy\ge \epsilon_0^2\}.$$
By
Proposition 2.4, we have that $u\in
W^{2,p}_{\hbox{loc}}(\Omega,S^k)$ for some $p>2$. In particular,
we have
$$\Sigma\subset\Sigma_p
=\{x\in\Omega | \ \liminf_{r\rightarrow 0}
r^{2p-n}\int_{B_r(x)}|\laplacian u|^p\,dy\ge \epsilon_1\}$$ for
some $\epsilon_1>0$. Therefore it is well-known (cf \cite{G}) that
${\Cal H}^{n-2p}(\Sigma)=0$. \qed \enddemo

Now we want to prove that a sequence  of weakly convergent
minimizing biharmonic maps is also  strongly convergent. More
precisely, we have
\proclaim {Lemma 3.3} For $n\ge 5$ and $k\ge
4$, let $\{u_i\}\subset W^{2,2}(\Omega,S^k)$ be a sequence of
minimizing biharmonic maps such that $u_i$ converges weakly in
$W^{2,2}$ to a map $u\in W^{2,2}(\Omega, S^k)$. Then $u_i$
converges to $u$ strongly in $W^{2,2}_{\hbox{loc}}(\Omega, S^k)$.
Moreover, $u\in W^{2,2}(\Omega,S^k)$ is also a minimizing
biharmonic map.
\endproclaim
\demo {Proof} The idea is similar to that of \cite{HL}. First, it
follows from Proposition 2.4 that there exists an $p>2$ such that
for any compact subset $K\subset\subset\Omega$
$$\sup_{i\ge 1}\|u_i\|_{W^{2,p}(K)}\le C(p,K)<\infty. \tag 3.8$$
Therefore, by the Rellich's compactness theorem, we may assume that
$u_i\rightarrow u$ strongly in $W^{1,4}_{\hbox{loc}}(\Omega,S^k)$.
By localization, it suffices to show that $u$ minimizes $\Bbb H$
on $B_R$ and $u_i\rightarrow u$ strongly in $W^{2,2}(B_R,S^k)$ for
any ball $B_{2R}\subset \Omega$.

Let $v\in W^{2,2}(B_{2R},S^k)$ be any map such that $v=u$ in
$B_{2R}\setminus B_R$. For any small $\delta>0$, let
$\eta_\delta\in C^\infty_0(B_{(1+3\delta)R})$ be such that $0\leq
\eta_\delta\leq 1$, $\eta_\delta\equiv 1$ in $B_{(1+2\delta)R}$,
$|\nabla\eta_\delta|\le {2\over \delta R}$, and
$|\nabla^2\eta_\delta|\le {4\over (\delta R)^2}$.

Consider $v_i(x)=\eta_\delta(x) v(x)+(1-\eta_\delta(x))u_i(x):
A_\delta\equiv B_{(1+3\delta)R} \setminus B_{(1+\delta)R}\to
\R^{k+1}$. Then it is easy to see $v_i\in W^{2,2}(A_\delta,
\R^{k+1})$ satisfies the condition of Lemma 2.1 on $A_\delta$.
Therefore Lemma 2.1 implies that $w_i\in W^{2,2}_{v_i}(A_{\delta},
S^k)$ such that
$$\int_{A_{\delta}} |\laplacian w_i|^2\,dx\leq
C\int_{ A_{\delta}}(| \laplacian v_i|^2+|\nabla v_i|^4)\,dx.$$
Let
$$\overline{w}_i(x)= \cases  w_i(x) &\text { for } x\in
A_{\delta}\\
v(x) &\text { for } x \in B_{(1+\delta)R}.\endcases$$
Then
${\overline {w}}_i\in W^{2,2}_{u_i}(B_{(1+3\delta)R}, S^k)$ so
that the $\Bbb H$-minimality of $u_i$ implies
$$\split\int_{B_{(1+3\delta)R}}|\laplacian u_i|^2\,dx &\leq \int_{B_{(1+3\delta)R}}
|\laplacian {\overline w}_i|^2\,dx\\
&= \int_{B_{(1+\delta)R}}|\laplacian v|^2\,dx
+\int_{A_{\delta}}|\laplacian w_i|^2\,dx\\
&\le \int_{B_{(1+\delta)R}}|\laplacian
v|^2\,dx+C\int_{A_\delta}(|\laplacian v_i|^2+ |\nabla
v_i|^4)\,dx.\endsplit \tag 3.9
$$
Direct calculations imply
$$\int_{A_{\delta}} |\nabla v_i|^4\,dx \leq
C[(\delta
R)^{-4}\int_{A_{\delta}}|u_i-v|^4\,dx+\int_{A_{\delta}}(|\nabla(u_i-v)|^4
+|\nabla v|^4)\,dx]$$

This, combined with the fact that $v=u$ on $A_\delta$ and
$u_i\rightarrow u$ strongly in $W^{1,4}$, implies
$$
\lim_{i\rightarrow \infty}\int_{A_\delta}|\nabla v_i|^4\,dx
=\int_{A_\delta}|\nabla v|^4\,dx=o(\delta),\tag 3.10
$$
where $\lim_{\delta\rightarrow 0}o(\delta)=0$. We also have
$$\split \int_{A_\delta}|\laplacian v_i|^2\,dx
&\leq C[(\delta R)^{-4}\int_{A_\delta}|u_i-v|^2\,dx
+(\delta R)^{-2}\int_{A_\delta}|\nabla(u_i-v)|^2\,dx]\\
&+C\int_{A_\delta}(|\laplacian u_i|^2+|\laplacian
v|^2)\,dx.\endsplit \tag 3.11
$$
Therefore, by (3.8), we have
$$
\split\lim_{i\rightarrow \infty}\int_{A_\delta}|\laplacian v_i|^2\,dx
&=\lim_{i\rightarrow\infty}\int_{A_\delta}|\laplacian u_i|^2\,dx
+\int_{A_\delta}|\laplacian v|^2\,dx\\
&\le \lim_{i\rightarrow \infty}(\int_{A_\delta}|\laplacian
u_i|^p)^{2\over p} |A_\delta|^{1-2\over p}
+\int_{A_\delta}|\laplacian v|^2\,dx\\
&=o(\delta).\endsplit \tag 3.12
$$
Putting (3.10), (3.11), and (3.12) into (3.9), we obtain
$$
\lim_{i\rightarrow\infty}\int_{B_{(1+3\delta)R}}|\laplacian
u_i|^2\,dx \le \int_{B_{(1+\delta)R}}|\laplacian
v|^2\,dx+o(\delta).\tag 3.13
$$
Letting $v\equiv u$ and $\delta\rightarrow 0$, (3.13) implies
$u_i\rightarrow u$ strongly in $W^{2,2}(B_R)$. Moreover, by the
lower semicontinuity, (3.13) also implies
$$\int_{B_{(1+3\delta)R}}|\laplacian u|^2\,dx\le
\int_{B_{(1+\delta)R}}|\laplacian v|^2\,dx+o(\delta)$$ this
clearly implies the $\Bbb H$-minimality of $u$ on $B_R$. The proof
is complete.   \qed
\enddemo

In order to give a proof of theorem A, we also need to recall the
following monotonicity inequality, which was established in
\cite{CWY} for stationary biharmonic maps.

 \proclaim{Lemma 3.4}
{For $n\ge 5$ and a compact Riemannian submanifold $N\subset
R^{k+1}$ without boundary. Suppose that $u\in W^{2,2}(\Omega, N)$
is a stationary biharmonic map. Then we have, for any $x\in
\Omega$ and $0<\rho\le r<\hbox{dist}(x,
\partial\Omega)$,
$$
\split & r^{4-n}\int_{B_r(x)}|\laplacian u|^2\,dy+r^{3-n}\int_{\partial B_r(x)}
[4|\nabla u|^2-4|{\partial u\over\partial r}|^2
+r{\partial\over\partial r}(|\nabla u|^2)]\,dy \\
&=\rho^{4-n}\int_{B_\rho(x)}|\laplacian u|^2\,dy
+\rho^{3-n}\int_{\partial B_\rho(x)} [4|\nabla u|^2-4|{\partial
u\over\partial \rho}|^2 +\rho{\partial\over\partial \rho}(|\nabla
u|^2
)]\,dy\\
&+4\int_{B_r(x)\setminus B_\rho(x)} ({|\nabla((y-x)\cdot\nabla
u)|^2\over |y-x|^{n-2}} +(n-2){|(y-x)\cdot\nabla u|^2\over
|y-x|^n})\,dy. \endsplit \tag 3.14
$$
\endproclaim

Now we complete a proof of theorem A.
\demo {Proof of Theorem A}
First it follows from Lemma 3.1 that the singular set ${\Cal
S}(u)$ is defined by
$${\Cal S}(u)=\{x\in\Omega: \Theta^{n-4}(u,x)\equiv
\liminf_{r\rightarrow 0}r^{4-n}\int_{B_r(x)}|\laplacian u|^2\,dy
\ge \epsilon_0^2\}.$$

It follows from Lemma 2.3 and (2.12) that there exists a $C=C(n,k)>0$ such that
$$\Theta^{n-4}(u,x)\leq \limsup_{r\rightarrow 0}r^{4-n}\int_{B_r(x)}|\laplacian u|^2\,dy
\le C, \ \forall x\in\Omega. \tag 3.15$$

Since minimizing biharmonic maps are stationary biharmonic maps,
Lemma 3.4 implies that for any $x\in\Omega$,
$$\sigma^u(x,r):=r^{4-n}\int_{B_r(x)}|\laplacian u|^2+r^{3-n}\int_{\partial B_r(x)}
[4|\nabla u|^2-4|{\partial u\over\partial
r}|^2+r{\partial\over\partial r}(|\nabla u|^2)]$$ is
monotonicially nondecreasing with respect to $r>0$ so that
$$\sigma^u(x)\equiv \lim_{r\rightarrow 0}\sigma^u(x,r)$$
exists for any $x\in\Omega$. It is easy to see
$\sigma^u(x)<+\infty$. To see $\sigma^u(x)>-\infty$, let $r>0$ be
a good slice, i.e.,
$$\split &|r^{3-n}\int_{\partial B_r(x)}r{\partial\over\partial r}(|\nabla u|^2)|\cr
&\le 2^n((2r)^{4-n}\int_{B_{2r}(x)}|\nabla^2 u|^2)^{1\over
2}((2r)^{2-n}\int_{ B_{2r}(x)}|\nabla u|^2)^{1\over 2}\le
C\endsplit$$ so that $\sigma^u(x,r)\ge -C$. Therefore, by choosing
good slices $r\downarrow 0$, we have $\sigma^u(x)\ge -C>-\infty$.

Next, we have

 \noindent{Claim:} {\it For any $x_0\in {\Cal S}(u)$
and $r_i\rightarrow 0$ there exists a minimizing biharmonic map
$\phi\in W^{2,2}_{\hbox{loc}}(\R^n, S^k)$ of homogeneous of degree
zero (i.e. $\phi(x)=\phi({x\over |x|})$) such that after passing
to subsequences $u_i(x)\equiv u(x_0+r_ix)$ converges to a
aminimizing biharmonic map $\phi$ strongly in
$W^{2,2}_{\hbox{loc}}(\R^n,S^k)$.}

To show this claim, it follows from (3.15) that for any $R>0$,
$\{u_i\}\subset W^{2,2}(B_R,S^k)$ is a bounded sequence of
minimizing biharmonic maps. Therefore, it follows from Lemma 3.3
that there exist a minimizing biharmonic map $\phi\in
W^{2,2}(B_R,S^k)$ such that $u_i\rightarrow \phi$ strongly in
$W^{2,2}(B_R,S^k)$. To see $\phi$ is homogeneous of degree zero,
note that for any $0<R_1<R_2\le R$
$$\split\sigma^{u_i}(0,R_2)-\sigma^{u_i}(0,R_1)
&=\sigma^u(x_0, R_2r_i)-\sigma^u(x_0,R_1r_i)\\
&\rightarrow \sigma^u(x_0)-\sigma^u(x_0)=0, \hbox{ as
}i\rightarrow\infty.\endsplit$$

This, combined with (3.14) and the lower semicontinuity, implies

$$\split &\int_{B_{R_2}\setminus B_{R_1}}({|\nabla(x\cdot\nabla \phi)|^2\over |x|^{n-2}}
+(n-2){|x\cdot\nabla \phi|^2\over |x|^n})\,dy\\
&\le \lim_{i\rightarrow \infty}\int_{B_{R_2}\setminus
B_{R_1}}({|\nabla(x\cdot\nabla u_i)|^2\over
|x|^{n-2}}+(n-2){|x\cdot\nabla u_i|^2\over
|x|^n})\,dy=0.\endsplit
$$
Therefore ${\partial\phi\over\partial r}=0$ for a.e. $x\in
B_{R_2}\setminus B_{R_1}$, which yields $\phi$ is of homogeneous
of degree zero.

With the help of Lemma 3.1, 3.3, 3.4, and the above claim, the
dimension estimation of ${\Cal S}(u)$ can be proved by a
refinement of the dimension reduction argument of Federer.   For
details, we refer to \cite {S; Chapter 3}. \qed
\enddemo

\head {\bf 4. Proof of Theorem B}\endhead

This section is devoted to the proof of Theorem B. One of the
crucial parts is to establish the sequentially lower
semicontinuity of $\Bbb H_\lambda(\cdot)$ in $W^{2,2}_\phi(\Omega,
S^4)$. Throughout this section we assume that $n=5$,
$\Omega\subset\R^5$, and $\phi\in W^{2,2}(\Omega,S^4)\cap
C^\infty(\partial\Omega,S^4)$.

Let's first recall the wedge product in $\R^5$. For the standard
orthonormal base $\{\bold e_i\}_{i=1}^5$ of $\R^5$, the wedge
product of four vectors $a$, $b$, $c$, $d$ in $\R^5$, $a\wedge
b\wedge c\wedge d\in \R^5$, is given by
$$
(a\wedge b\wedge c\wedge d)_i=\text {det } (\bold e_i, a,b,c,d), 1\leq i\leq 5.
$$

Now we have \proclaim {Lemma 4.1} For every $\lambda\in (0,1]$,
$\Bbb H_{\lambda}(\cdot)$ is s.l.s.c. in $W^{2,2}_\phi(\Omega,
S^4)$ for the weak $W^{2,2}$ topology.
\endproclaim
\demo {Proof} Since the supremum of sequentially lower
semicontinuous functions is still a sequentially lower
semicontinuous function, it suffices to prove that for any fixed
$\xi :\Omega\to \R$ with $\|\xi\|_{L^{\infty}}\leq 1$ the
functional
$$
\Bbb H_{\lambda, \xi} (u)=\int_{\Omega}|\laplacian
u|^2\,dx+16 \lambda \int_{\Omega} D(u)\cdot\nabla \xi\,dx
$$ is sequentially lower semicontinuous in $W^{2,2}_\phi(\Omega,S^4)$ for
the weak $W^{2,2}$ topology.

Let $\{u^n\}\subset W^{2,2}_\phi(\Omega,S^4)$ converge to $u\in
W^{2,2}_\phi(\Omega,S^4)$ weakly in $W^{2,2}(\Omega,S^4)\cap
W^{1,4}(\Omega,S^4)$, and strongly in $W^{1,2}(\Omega,S^4)$. Set
$v^n=u^n-u\in W^{2,2}(\Omega,\R^5)$. Then we have
$$\int_{\Omega} |\laplacian u^n|^2\,dx=\int_{\Omega} |\laplacian
v^n|^2\,dx +\int_{\Omega} |\laplacian u|^2\,dx+o(1),\tag 4.1$$
where $o(1)$ is such that $\lim_{n\to \infty}o(1)=0$.

Now we claim
$$\int_{\Omega} |\nabla v^n|^4\,dx\leq \int_{\Omega}
|\laplacian v^n|^2\,dx +o(1). \tag 4.2
$$

To show (4.2), observe that since $|u^n|=1$ and $|u|=1$, we have
$$|v^n|^2=-2v^n\cdot u$$
this implies, by taking two derivatives,
$$|\nabla v^n|^2=-(\laplacian v^n\cdot v^n+\laplacian v^n\cdot u
+ v^n\cdot\laplacian u +2\nabla v^n\cdot\nabla u).$$

Since we have, for a.e. $x\in\Omega$,
$$(\laplacian u \cdot v^n)\rightarrow 0, \ (\laplacian u \cdot v^n)u^n\rightarrow 0,
\ (\laplacian u \cdot v^n)\nabla u\rightarrow 0$$ and

$$\max\{|(\laplacian u \cdot v^n)|, \ |(\laplacian u \cdot v^n)u^n|\}
\leq 2|\laplacian u|\in L^2(\Omega), \ |(\laplacian u \cdot
v^n)\nabla u|\le |\laplacian u||\nabla u|\in L^{4\over
3}(\Omega).$$
The Lebegues Dominated Convergence Theorem implies
$$\lim_{n\rightarrow\infty}\int_{\Omega}
(|\laplacian u\cdot v^n|^2+|(\laplacian u\cdot v^n)u^n|^2+
|(\laplacian u\cdot v^n)\nabla u|^{4\over 3})\,dx=0$$ so that we
have
$$\lim_{n\rightarrow\infty}\int_{\Omega}
(\laplacian u\cdot v^n)u^n\cdot \laplacian v^n\,dx=0, \
\lim_{n\rightarrow \infty}\int_{\Omega} (\laplacian u\cdot
v^n)\nabla u\cdot \nabla v^n\,dx=0.$$
Now we need to show
$$\lim_{n\rightarrow\infty}\int_{\Omega}|\nabla u\cdot\nabla v^n|^2\,dx=0.\tag 4.3$$
Assume that (4.3) is true for the moment. Then it is easy to see
$(\nabla v^n\cdot u^n)(\nabla u\cdot\nabla v^n)\rightarrow 0$ in
$L^1(\Omega)$ so that we have
$$\int_\Omega |\nabla v^n|^4=\int_\Omega |u^n\cdot\laplacian v^n|^2\,dx
+o(1)\leq\int_\Omega|\laplacian v^n|^2\,dx+o(1)$$ this clearly
implies (4.2). To see (4.3), observe that $v^n\in
W^{2,2}_0(\Omega,\R^5)$. Therefore we have by integration by
parts,
$$\split \int_{\Omega}|\nabla u\cdot \nabla v^n|^2\,dx
&=\int_\Omega (\nabla u\cdot\nabla v^n)(\nabla u\cdot\nabla v^n)\,dx\\
&=-\int_\Omega \nabla\cdot((\nabla u\cdot\nabla v^n)\nabla u)\cdot v^n\,dx\\
&=-2\int_\Omega (\laplacian u\cdot\nabla v^n)\nabla u\cdot v^n\,dx
-\int_\Omega (\nabla u\cdot\laplacian v^n)(\nabla u\cdot v^n)\,dx\\
&\le C\int_\Omega (|\nabla v^n||\nabla u||\laplacian u|+|\nabla
u|^2|\laplacian v^n|)\,dx \rightarrow 0, \hbox{ as }
n\rightarrow\infty.\endsplit$$
This gives (4.3).

Now we write
$$\int_{\Omega} D(u^n)\cdot\nabla\xi =A^n+B^n+C^n,$$
where

$$A^n=\int_{\Omega} u^n\cdot [(\frac {\partial u}{\partial x_2}
\wedge \cdots\wedge \frac {\partial u}{\partial x_5}) \frac
{\partial \xi}{\partial x_1}+\cdots+(\frac {\partial u}{\partial
x_1}\wedge \cdots\wedge \frac {\partial u}{\partial x_4})\frac
{\partial \xi}{\partial x_5}]\,dx,$$
$$B^n=\int_{\Omega} [D(u^n)\cdot\nabla\xi\,dx- A^n-C^n$$
and
$$C^n=\int_{\Omega} V^n\cdot\nabla\xi\,dx,
$$
where $$V^n=(\text {det }(u^n,\frac {\partial v^n}{\partial x_2},
\cdots,\frac {\partial v^n}{\partial x_5}), \cdots, \text { det }
(u^n,\frac {\partial v^n}{\partial x_1}, \cdots, \frac {\partial
v^n}{\partial x_4})).$$

Since $u^n\rightarrow u$ weak$^*$ in $L^\infty(\Omega)$, we have
$$A^n\to\int_{\Omega}D(u)\cdot \nabla \xi\,dx, \hbox{ as }n\rightarrow\infty.$$

To estimate $B^n$, we observe that direct calculations imply
$$\split
&\quad \frac {\partial u^n}{\partial x_2}\wedge \frac {\partial
u^n}{\partial x_3}\wedge\frac {\partial u^n}{\partial x_4}
\wedge\frac {\partial u^n}{\partial x_5}-\frac {\partial
v^n}{\partial x_2}\wedge \frac {\partial v^n}{\partial
x_3}\wedge\frac {\partial v^n}{\partial
x_4} \wedge\frac {\partial v^n}{\partial x_5}\\
&=\frac {\partial u}{\partial x_2}\wedge \frac {\partial
u}{\partial x_3}\wedge\frac {\partial u}{\partial x_4} \wedge\frac
{\partial u}{\partial x_5}+\frac {\partial u}{\partial x_2}\wedge
\frac {\partial u}{\partial x_3}\wedge\frac {\partial u}{\partial
x_4} \wedge\frac {\partial v^n}{\partial x_5} +\frac {\partial
u}{\partial x_2}\wedge \frac {\partial u}{\partial x_3}\wedge\frac
{\partial v^n}{\partial x_4} \wedge\frac {\partial v^n}{\partial
x_5}\\
&+\frac {\partial u}{\partial x_2}\wedge \frac {\partial
u}{\partial x_3}\wedge\frac {\partial v^n}{\partial x_4}
\wedge\frac {\partial u}{\partial x_5}+\frac {\partial u}{\partial
x_2}\wedge \frac {\partial v^n}{\partial x_3}\wedge\frac {\partial
u}{\partial x_4} \wedge\frac {\partial u}{\partial x_5}+\frac
{\partial u}{\partial x_2}\wedge \frac {\partial v^n}{\partial
x_3}\wedge\frac {\partial u}{\partial x_4} \wedge\frac {\partial
v^n}{\partial x_5}\\
&+\frac {\partial v^n}{\partial x_2}\wedge \frac {\partial
v^n}{\partial x_3}\wedge\frac {\partial v^n}{\partial x_4}
\wedge\frac {\partial u }{\partial x_5}+\frac {\partial
u}{\partial x_2}\wedge \frac {\partial v^n}{\partial
x_3}\wedge\frac {\partial v^n}{\partial x_4} \wedge\frac {\partial
v^n}{\partial x_5}+\frac {\partial v^n}{\partial x_2}\wedge \frac
{\partial u}{\partial x_3}\wedge\frac {\partial v^n}{\partial x_4}
\wedge\frac {\partial v^n}{\partial x_5}\\
&+\frac {\partial v^n}{\partial x_2}\wedge \frac {\partial
u}{\partial x_3}\wedge\frac {\partial v^n}{\partial x_4}
\wedge\frac {\partial u}{\partial x_5}+\frac {\partial
v^n}{\partial x_2}\wedge \frac {\partial u}{\partial
x_3}\wedge\frac {\partial u}{\partial x_4} \wedge\frac {\partial
v^n}{\partial x_5}+\frac {\partial v^n}{\partial x_2}\wedge \frac
{\partial u}{\partial x_3}\wedge\frac {\partial u}{\partial x_4}
\wedge\frac {\partial u}{\partial x_5}\\
&+\frac {\partial v^n}{\partial x_2}\wedge \frac {\partial
v^n}{\partial x_3}\wedge\frac {\partial u}{\partial x_4}
\wedge\frac {\partial v^n}{\partial x_5}+\frac {\partial
v^n}{\partial x_2}\wedge \frac {\partial v^n}{\partial
x_3}\wedge\frac {\partial u}{\partial x_4} \wedge\frac {\partial
u}{\partial x_5}+\frac {\partial v^n}{\partial x_2}\wedge \frac
{\partial v^n}{\partial x_3}\wedge\frac {\partial v^n}{\partial
x_4} \wedge\frac {\partial u }{\partial x_5}
\endsplit $$

This implies
$$|B^n|\leq C\int_\Omega (|\nabla u|^3 |\nabla v^n|+|\nabla u|^2|\nabla v^n|^2
+|\nabla u||\nabla v^n|^3)\,dx.$$

Since $\nabla v^n\rightarrow 0$ a.e. and weakly in $L^4(\Omega)$,
we can conclude that $\nabla v^n\rightarrow 0$ strongly in
$L^q(\Omega)$ for any $1\le q<4$. Therefore we see $B^n\rightarrow
0$ as $n\rightarrow\infty$.

Now we need to estimate $C^n$. Since $|u(x)|=1$ for a.e.
$x\in\Omega$, there exists a rotation $R\in SO(4)$ such that
$R(u(x))=(0,0,0,0,1)$. Moreover, for any vector $p_1,
p_2,p_3,p_4,p_5)\in \R^5$, since
$$\text {det}(u(x), p_{i_1}, \cdots, p_{i_4})
=\text {det}(R(u(x)), R(p_{i_1}), \cdots, R(p_{i_4})), \forall
1\le i_1<i_2<i_3<i_4\le 5.$$ We may assume that $u(x)=(0,0,0,0,1)$
and write $R(p_i)=(a_i,b_i,c_i, d_i, e_i)$ for $1\le i\le 5$. Then
we have
$$V=(\text {det}(u(x), p_2,\cdots,
p_5),\cdots, \text {det}(u(x), p_1, \cdots, p_4)) =a\wedge b\wedge
c\wedge d,$$ where $a=(a_1,\cdots, a_5), b=(b_1,\cdots,
b_5),\cdots, d=(d_1,\cdots, d_5)$. Therefore we have
$$|V|=|a\wedge b\wedge c\wedge d|\le {1\over 2^4}(|a|^2+|b|^2+|c|^2+|d|^2)^2. \tag 4.4$$
Applying (4.4) with $p_i={\partial v^n\over\partial x_i}$ for
$1\le i\le 5$, we obtain
$$|V^n|(x)\le {1\over 16}|\nabla v^n|^4(x), \hbox{ for a.e. }x\in\Omega. \tag 4.5$$
This, combined with (4.1) and (4.2), implies
$$\liminf_{n\to\infty}[\int_{\Omega} |\laplacian u^n|^2dx+16\lambda
\int_\Omega D(u^n)\cdot\nabla\xi \,dx]\geq \int_{\Omega}
|\laplacian u|^2dx+16\lambda \int_\Omega D(u)\cdot\nabla\xi
\,dx.$$ This completes the proof of Lemma 4.1.       \qed
\enddemo

As a direct consequence, we have
\proclaim{Corollary 4.2} For any
$\lambda\in (0,1]$, there exists a $u_\lambda\in
W^{2,2}_\phi(\Omega,S^4)$ which minimizes $\Bbb H_\lambda(\cdot)$
over $W^{2,2}_\phi(\Omega,S^4)$.
\endproclaim
\demo {Proof} Since $L(u)\ge 0$ for any $u\in
W^{2,2}_\phi(\Omega,S^4)$, it is easy to see that any minimizing
sequence $\{u_i\}$ of $\Bbb H_\lambda(\cdot)$ over
$W^{2,2}_\phi(\Omega,S^4)$ is a bounded sequence in
$W^{2,2}(\Omega)$. Therefore we may assume that $u_i$ converges to
$u_\lambda$ weakly in $W^{2,2}(\Omega)$. By Lemma 4.1, we have
that $u_\lambda$ is a minimizer for $\Bbb H_\lambda$ over
$W^{2,2}_\phi(\Omega,S^4)$.  \qed \enddemo

Now we have
 \proclaim{Lemma 4.3} For any $\lambda\in (0,1)$, if
$u_\lambda\in W^{2,2}_\phi(\Omega,S^4)$ is a minimizer for $\Bbb
H_\lambda(\cdot)$. Then $u_\lambda$ is a $Q$-minimizing biharmonic
map, with $Q={1+\lambda\over 1-\lambda}$.
\endproclaim
\demo {Proof} For simplicity, we abbreviate $u_{\lambda}$ to $u$.
Let $w\in W^{2,2}_\phi(\Omega, S^4)$. Then, similar to (4.5), we
have
$$|D(w)|(x)\leq \frac 1{16} |\D w|^4(x), \hbox{ for a.e. }x\in\Omega,$$
and
$$\split |L(w)-L(u)|&\leq \int_\Omega (|D(w)|+|D(u)|)\,dx\\
&\leq \frac 1 {16\sigma_4}\int_{\Omega}(|\nabla w|^4 +|\nabla
u|^4)\,dx\leq \frac 1 {16\sigma_4} [\Bbb H(w)+\Bbb H(u)],\endsplit
\tag 4.6$$ where we have used the fact that $|\nabla w|^2\le
|\laplacian w|^2$ and $|\nabla u|^2\le |\laplacian u|^2$ for a.e.
$x\in\Omega$.

Since $u$ minimizes $\Bbb H_\lambda$, we then have
$$\Bbb H(u)\leq \Bbb H(w)+ 16\lambda \sigma_4 (L(w)-L(u))
\le \Bbb H(w)+\lambda (\Bbb H(w)+\Bbb H(u)).$$

This implies
$$\Bbb H(u)\leq \frac
{1+\lambda} {1-\lambda } \Bbb H (w), \forall w\in
W^{2,2}_{\phi}(\Omega , S^4). \tag 4.7$$

Now we need to show that $u_\lambda $ is a biharmonic map. To see
it, let $\eta\in C_0^\infty(\Omega, \R^5)$, $t\in [0,1)$, and
denote $u_\lambda^t(x)={u(x)+t\eta(x)\over |u(x)+t\eta(x)|}$ for
$x\in\Omega$. Then we have
$${d\over dt}|_{t=0}(\Bbb H(u_\lambda^t)+16\lambda \sigma_4 L(u_\lambda^t))=0. \tag 4.8$$
Therefore $u_\lambda$ is a biharmonic map, if we can show
$${d\over dt}|_{t=0}  L(u_\lambda^t))=0. \tag 4.9$$

In order to prove (4.9), we need the following Lemmas.
\proclaim{Lemma 4.4} For any $u,v\in W^{2,2}_\phi(\Omega,S^4)$, we
have the following inequality
$$|L(u)-L(v)|\le C\|\nabla(u-v)\|_{L^4(\Omega)}
(\|\nabla u\|_{L^4(\Omega)}^3+\|\nabla v\|_{L^4(\Omega)}^3). \tag
4.10$$
\endproclaim

\demo{Proof} By the definition of $L$, we see
$$|L(u,u_0)-L(v,u_0)|\leq L(u,v)=\frac 1{\sigma_4}
\sup_{\xi :\Omega\to \R; |\nabla \xi|\leq 1} \int_{\Omega}
(D(u)-D(v))\cdot \nabla\xi\,dx\tag 4.11
$$

For any $\xi:\Omega\to R$ with $|\nabla \xi|\le 1$, we write
$$ \int_{\Omega} (D(u)-D(v))\cdot\nabla\xi\,dx =I+II+III+IV
+V,\tag 4.12
$$
where
$$\split
I=\int_{\Omega}[ \text {det} &\left (u-v,\frac {\partial
u}{\partial x_2}, \frac {\partial u}{\partial x_3},\frac {\partial
u}{\partial x_4},   \frac {\partial u}{\partial x_5} \right )\frac
{\partial \xi}{\partial x_1}+\cdots \\
&+\text {det} \left (u-v, \frac {\partial u}{\partial x_1},
 \frac {\partial u}{\partial x_2},\frac {\partial u}{\partial x_3},
\frac {\partial u}{\partial x_4}\right )\frac {\partial
\xi}{\partial x_5}]\,dx,\endsplit
$$
$$\split
II=\int_{\Omega}[\text {det} &\left (v,\frac {\partial
(u-v)}{\partial x_2}, \frac {\partial u}{\partial x_3},\frac
{\partial u}{\partial x_4},
 \frac {\partial u}{\partial x_5} \right )\frac {\partial \xi}{\partial
x_1}+\cdots \\
&+\text {det} \left (v,\frac {\partial (u-v)}{\partial x_1}, \frac
{\partial u}{\partial x_2},\frac {\partial u}{\partial x_3}, \frac
{\partial u}{\partial x_4}\right )\frac {\partial \xi}{\partial
x_5}]\,dx,
\endsplit
$$
$$\split
III=\int_{\Omega} [\text {det} &\left (v,\frac {\partial
v}{\partial x_2}, \frac {\partial (u-v)}{\partial x_3},\frac
{\partial u}{\partial x_4}, \frac {\partial u}{\partial x_5}
\right )\frac {\partial \xi}{\partial
x_1}+\cdots \\
&+\text {det} \left (v,\frac {\partial v}{\partial x_1},
 \frac {\partial (u-v)}{\partial x_2},\frac {\partial u}{\partial x_3},
 \frac {\partial u}{\partial x_4}\right )\frac {\partial
\xi}{\partial x_5}]\,dx,\endsplit
$$
$$\split
IV=\int_{\Omega}[\text {det} &\left (v,\frac {\partial v}{\partial
x_2}, \frac {\partial v}{\partial x_3},\frac {\partial
(u-v)}{\partial x_4},
 \frac
{\partial u}{\partial x_5} \right )\frac {\partial \xi}{\partial
x_1}+\cdots \\
&+\text {det} \left (v,\frac {\partial v}{\partial x_1},
 \frac {\partial v}{\partial x_2},\frac {\partial (u-v)}{\partial x_3},
 \frac {\partial u}{\partial x_4}\right )\frac {\partial
\xi}{\partial x_5}]\,dx,
\endsplit
$$
$$\split V=\int_{\Omega} [\text {det} &\left (v,\frac {\partial v}{\partial x_2},
\frac {\partial v}{\partial x_3},\frac {\partial v}{\partial x_4},
 \frac
{\partial (u-v)}{\partial x_5} \right )\frac {\partial
\xi}{\partial x_1}+\cdots \\
&+\text {det} \left (v,\frac {\partial u}{\partial x_1},
 \frac {\partial v}{\partial x_2},\frac {\partial v}{\partial x_3},
 \frac {\partial (u-v)}{\partial x_4}\right )\frac {\partial
\xi}{\partial x_5}]\,dx.
\endsplit $$

It follows from H\"older's inequality that
$$|II|\leq C\int_{\Omega} |\nabla (u-v)||\nabla u|^3\,dx,$$

$$\split|III|&\leq C\int_{\Omega} |\nabla (u-v)||\nabla
u|^2|\nabla v|\,dx\\
&\leq C(\int_{\Omega} |\nabla (u-v)|^4\,dx )^{1/4} (\int_{\Omega}
(|\nabla u|^4+|\nabla v|^4)\,dx)^{1/4},\endsplit
$$

$$\split |IV|&\leq C\int_{\Omega} |\nabla (u-v)||\nabla u||\nabla
v|^2\,dx\\
&\leq C(\int_{\Omega} |\nabla (u-v)|^4\,dx )^{1/4} (\int_{\Omega}
(|\nabla u|^4+|\nabla v|^4)\,dx)^{1/4},\endsplit
$$
and
$$|V|\leq C\int_{\Omega} |\nabla (u-v)||\nabla v|^3 \,dx\leq
C(\int_{\Omega} |\nabla (u-v)|^4\,dx )^{1/4} (\int_{\Omega}
|\nabla v|^4\,dx)^{1/4}.$$

In order to estimate $I$, we observe that
$$\split
&\quad 4  \left (\frac {\partial u}{\partial x_2}\wedge  \frac
{\partial u}{\partial x_3}\wedge \frac {\partial u}{\partial x_4}
\wedge \frac {\partial u}{\partial x_5}\right )\\
 =&\left (   u \wedge \frac {\partial u}{\partial x_3}
 \wedge \frac {\partial u}{\partial x_4}
 \wedge \frac
{\partial u}{\partial x_5}\right )_{x_2}+  \left (\frac {\partial
u}{\partial x_2}\wedge u\wedge \frac {\partial u}{\partial x_4}
\wedge
\frac {\partial u}{\partial x_5}\right )_{x_3}\\
 &+  \left (\frac {\partial
u}{\partial x_2}\wedge \frac {\partial u}{\partial x_3}\wedge
 u \wedge \frac {\partial
u}{\partial x_5}\right )_{x_4}+  \left (\frac {\partial
u}{\partial x_2}\wedge \frac {\partial u}{\partial x_3}\wedge
\frac {\partial u}{\partial x_4} \wedge  u\right )_{x_5}
\endsplit$$
in the sense of distributions. Therefore, by integration by parts,
we have
$$\split
&\int_{\Omega}  (u-v)\cdot\frac {\partial u}{\partial x_2}\wedge
\frac {\partial u}{\partial x_3}\wedge\frac {\partial u}{\partial
x_4}\wedge  \frac {\partial u}{\partial x_5}  \,\frac {\partial
\xi}{\partial x_1}\\
&=\frac 14\int_{\Omega} [\frac {\partial (u-v)}{\partial x_2}\cdot
u\wedge  \frac {\partial u}{\partial x_3}\wedge \frac {\partial
u}{\partial x_4}\wedge
\frac {\partial u}{\partial x_5} +\cdots\\
&+\frac {\partial (u-v)}{\partial x_5}\cdot  \frac {\partial
u}{\partial x_2}\wedge \frac {\partial u}{\partial x_3}\wedge
\frac {\partial u}{\partial x_4}\wedge  u ]\frac {\partial
\xi}{\partial x_1}\,dx\\
&+ \frac 1 4\int_{\Omega}(u-v)\cdot [(u\wedge\frac {\partial
u}{\partial x_3}\wedge \frac {\partial u}{\partial x_4}\wedge
\frac {\partial u}{\partial x_5} )\frac {\partial^2 \xi}{\partial
x_1\partial
x_2}+\cdots\\
& + (\frac {\partial u}{\partial x_2}\wedge \frac {\partial
u}{\partial x_3}\wedge \frac {\partial u}{\partial x_4}\wedge u
)\frac {\partial^2 \xi}{\partial x_1\partial x_5}]\,dx.
\endsplit$$
By doing the same calculations to all other terms in $I$, we see
that the sum of all terms involving $\nabla^2\xi$ cancel each
other. Therefore we have
$$|I|\leq C\int_{\Omega} |\nabla (u-v) |\nabla u|^3\,dx\leq \left
(\int_{\Omega} |\nabla (u-v)|^4\,dx \right )^{1/4}\left
(\int_{\Omega}|\nabla u|^4\,dx\right )^{3/4}.$$ Putting all these
inequalities together, we obtain (4.10). \qed \enddemo

The next Lemma is concerning with the density of maps, which are
smooth away from finitely many singular points, in
$W^{2,2}_\phi(\Omega,S^4)$. The proof will be given in Section 5.
\proclaim{Lemma 4.5} Define
$$R_{\phi}^\infty=\{u\in W^{2,2}_\phi(\Omega,S^4):
u\in C^\infty(\bar\Omega\setminus \cup_{i=1}^l\{a_i\},S^4), \text{
where } l<\infty \text{ and } \cup_{i=1}^l
\{a_i\}\subset\Omega\}.$$
Then $R_{\phi}^\infty$ is dense in
$W^{2,2}_\phi(\Omega,S^4)$ for the $W^{2,2}$-topology.
\endproclaim

Now we return to the proof of Lemma 4.3. First we observe that for
any $v\in R_{\phi}^\phi$ we have
$$L({v+t\eta\over |v+t\eta|})=L(v), \hbox{ for sufficiently small } t\in [0,1)$$
since the singularity of $v^t={v+t\eta\over |v+t\eta|}\in
R_{\phi}^\infty$ is same as that of $v$ and $L(\cdot)$ is the
minimal connection of its singular points (\cite{BB} \cite{BCL}
\cite{BBC} \cite{GMS}).

For $u_\lambda$, it follows from Lemma 4.5 that there are
$\{u_n\}\subset R_\phi^\infty$ such that
$$\lim_{n\rightarrow \infty}\|u_n-u_\lambda\|_{W^{2,2}(\Omega)}=0.$$
Then, for sufficiently small $t\in [0,1)$, we also have
$$\lim_{n\rightarrow\infty}\|u_n^t-u_\lambda^t\|_{W^{2,2}(\Omega)}=0$$
where $u_n^t={u_n+t\eta\over |u_n+t\eta|}$. By Lemma 4.4, we have
$$\lim_{n\rightarrow\infty}L(u_n^t)=L(u_\lambda^t), \
\lim_{n\rightarrow\infty}L(u_n)=L(u_\lambda).$$

On the other hand, since $u_n\in R_\phi^\infty$, we have, for any
$t\in [0,1)$ sufficiently small,
$$L(u_n^t)=L(u_n).$$
Therefore we have $L(u_\lambda^t)=L(u)$ for any sufficiently small
$t\in [0,1)$. This finishes the proof of Lemma 4.3.    \qed
\enddemo

We completion of proof of Theorem B.

\demo {Proof of theorem B} Part (a) and (b) follow from Lemma 4.1
and Corollary 4.2. Since Lemma 4.3 implies that each $\Bbb
H_\lambda$-minimizer $u_\lambda$ is a $Q$-minimizing biharmonic
map with $Q={1+\lambda\over 1-\lambda}$, part (c) follows from
Corollary 3.2. \qed
\enddemo

\head {\bf 5. Appendix}\endhead

In this section, we provide two examples, a proof of Lemma 4.5, a
boundary partial regularity for $\Bbb H_\lambda$, and propose a
few open questions.

\proclaim {Proposition A1} For $n\ge 5$, $\Phi(x)={x\over |x|}:
B^n\to S^{n-1}$ is a unique minimizing biharmonic map in
$W^{2,2}_\Phi(B^5,S^4)$.
\endproclaim
\demo{Proof} Using the fact that $|u|^2=1$, we have
$$-u\cdot\laplacian u=|\nabla u|^2\tag 5.1
$$
Then
$$|\laplacian u|^2=|u\cdot\laplacian u|^2+|\laplacian u\cdot \tau
(u)|^2,
$$
where $\tau (u)=(\tau_1,...,\tau_{n-1})$ and $\{\tau_k(u)\}$ is
an orthonormal base of the tangent plane of $S^{n-1}$ at $u$. Since
$\Phi:B^n\to S^{n-1}$ is a weakly harmonic map,
$$\laplacian \Phi (x) \cdot  \tau (\Phi (x) )=0. $$
Now we recall that $\Phi :B^n\to S^{n-1}$ is a unique minimizing
$4$-harmonic map (cf. \cite {CG}, \cite {AL}, \cite {Ho}). Then
$$\split
\int_{\Omega} |\laplacian u|^2&=\int_{\Omega} |\nabla
u|^4+|\laplacian u\cdot\tau (u)|^2\,dx\\
&\geq \int_{\Omega} |\nabla \frac x{|x|}|^4\,dx = \int_{\Omega}
|\laplacian\frac  x{|x|}|^2\,dx
\endsplit$$
for all $u\in W^{2,2}_{\Phi}
(B^n; S^{n-1})$. This implies that $\Phi$ is a unique minimizing
biharmonic map. \qed  \enddemo

Now we give an example consisting a domain $\Omega\subset \R^5$
and $\phi:\partial\Omega\to S^4$ with $\deg (\phi)=0$ such that
the infimum of $\Bbb H$ in $W^{2,2}_\phi(\Omega,S^4)$ is less than
that in $C^\infty_\phi(\bar\Omega,S^4)$.

For a sufficiently large $L>0$, let $B_1^+((0',L))$ ( or
$B_1^-((0',-L))$ resp.) be the upper (or lower, resp.) half unit
ball centered at $(0',L)$ (or $(0',-L)$ resp.) in $\R^5$. Define

$$\Omega=B_1^+((0',L))\cup (B_1^4\times [-L, L])\cup B_1^-((0',-L))$$
where $B_1^4\subset \R^4$ is the unit ball centered at $0'\in \R^4$.

Let $\psi^+: \partial B_1^+((0', L))\cap\{x\in \R^5: x_5>L\}\to
S^4$ be a smooth map of degree one such that

$$\psi^+|_{\partial B_1^+((0',L))\cap\{x\in \R^5: x_5=L\}}=(0', 1),
\ {\partial\psi^+\over\partial x_5}|_{\partial
B_1^+((0',L))\cap\{x\in \R^5: x_5=L\}}=0.$$

Define $\psi:\partial\Omega\to S^4$ by

$$\psi(x',x_5)=\cases \psi^+(x',x_5), \  x_5\ge L\\
                 (0',1)  \ \ \ \ \ , \ x_5\in [-L,L] \\
                  \psi^+(x', -x_5), \ x_5\le -L. \endcases $$

Then $\psi\in C^0(\partial\Omega,S^4)\cap
W^{2,2}(\partial\Omega,S^4)$ has $\hbox{deg}(\phi)=0$. Motivated
by the gap phenomena discovered by Hardt-Lin \cite{HL} in the
context of harmonic maps, we have

\proclaim{Proposition A2} Under the above notations, we have the
following gap phenomena
$$\inf_{u\in W^{2,2}_\psi(\Omega,S^4)}\Bbb H(u)
<\inf_{v\in W^{2,2}_\psi(\Omega,S^4)\cap C^0(\bar\Omega,S^4)} \Bbb
H(v). \tag 5.5$$
\endproclaim

\demo{Proof} The idea is similar to that of \cite{HL1}. First,
observe that
$$\Psi(x)=\cases \psi({x-(0',L)\over |x-(0',L)|}), \ \ \ x\in B_1^+((0',L))\\
          \psi({x-(0',-L)\over |x-(0',-L)|}), \ \ x\in B_1^-((0',-L))\\
           (0',1), \ \ \ \ \ \ \ \ \ \ \ x\in B_1^4\times [-L,L].\endcases$$

Then it is not difficult to verify that $\Psi\in
W^{2,2}_\psi(\Omega,S^4)$. Moreover, direct calculations imply

$$\split\Bbb H(\Psi)
&=2\int_{B_1^+((0',1))}|\laplacian (\psi({x-(0',L)\over |x-(0',L)|}))|^2\,dx\\
&=2\int_0^1 \,dr \int_{\partial B_1^+(0',L)}|\laplacian \psi^+|^2
=C(\psi^+)\endsplit$$ is independent of $L$.

On the other hand, for any $v\in W^{2,2}_\psi(\Omega,S^4) \cap
C(\bar\Omega,S^4)$, since $|\laplacian v|(x)\ge |\nabla v|^2(x)$
for a.e. $x\in\Omega$, we have
$$\split &\Bbb H(v) \ge \int_{B_1^4\times [-L,L]}|\laplacian v|^2\,dx
\ge \int_{-L}^L\,dx_5 \int_{B_1^4}|\laplacian v|^2\,dx'\\
&\ge \int_{-L}^L \,dx_5\int_{B_1^4}|\nabla_{x'} v|^4(x',x_5)\,dx'
\ge 16 \int_{-L}^L\,dx_5\int_{B_1^4}|\hbox{det}(\nabla_{x'}v|(x',x_5)\,dx'\\
&\ge 16 \sigma_4\int_{-L}^L\,dx_5=32\sigma_4 L \endsplit \tag
5.6$$ where we have used the inequality (4.5) and
$$\int_{B_1^4}|\hbox{det}(\nabla_{x'}v)|(x',x_5)\,dx'\ge \sigma_4,
 \ \forall x_5\in (-L,L). \tag 5.7$$
(5.7) holds, since for any $x_5\in (-L,L)$ $v\in
C(\overline{B_1^+(0',L)\cup B_1^4\times [x_5,L]},S^4)$ hence
$v:\partial(B_1^+(0',L)\cup B_1^4\times [x_5,L])\to S^4$ has
degree zero.

In particular, we have that $v(\cdot,x_5): B_1^4\to S^4$ has
degree one for all $x_5\in (-L,L)$.

Therefore we establish (5.5), provided that $L>0$ is chosen to be
sufficiently large. \qed
\enddemo

Now,we complete the proof of Lemma 4.5.

\demo{Proof of Lemma 4.5} The idea is similar to that of
\cite{BZ}. Since $\phi\in C^\infty(\bar{\Omega}\setminus
\{x_i\}_{i=1}^k,S^4)$ for some $\{x_i\}_{i=1}^k\subset\Omega$, we
have that for any $u\in W^{2,2}_\phi(\Omega,S^4)$ there are
$\{u_n\}\subset C^\infty(\bar\Omega,\R^5)$ such that $u_n=\phi$,
$\nabla u_n=\nabla \phi$ on $\partial\Omega$, and {$u_n\rightarrow
u$ strongly in $W^{2,2}(\Omega,S^4)$. For any small $\epsilon>0$,
set
$$S^4_{1-\varepsilon }=\{ x\in \R^5 : |x|=1-\varepsilon \},
\quad S^4_{1+\varepsilon }=\{ x\in \R^5 : |x|=1+\varepsilon \}.$$

By the Sard's theorem, we have that
$$ F^-_{n,\varepsilon}=u^{-1}_n(S^4_{1-\varepsilon }),
\quad F^+_{n,\varepsilon}=u^{-1}_n(S^4_{1+\varepsilon })
$$
are two compact submanifolds of $\Omega$ of codimension one.
Moreover
$$V_{n,\varepsilon}^-=u_n^{-1}(\{ |y|\leq 1-\varepsilon \}),
\quad V_{n,\varepsilon}^+=u_n^{-1}(\{ |y|\geq 1+\varepsilon \})$$
are smooth domains inside $\Omega$ such that $\partial
V_{n,\epsilon}^{-}= F^-_{n,\varepsilon}$ and
 $\partial
V_{n,\epsilon}^+=F^+_{n,\varepsilon}$.

For any $a\in B^5_{1\over 2}$, define the projection maps $
\Pi_a^-:\R^5\to S^4_{1-\varepsilon}$  and $\Pi_a^+:\R^5\to
S^4_{1+\varepsilon}$ by
 $$\Pi_{a}^-(x)=\frac {x-a}{|x-a|} (1-\varepsilon ),\quad
\Pi_{a}^+(x)=\frac {x-a}{|x-a|} (1+\varepsilon ).
$$

By Lemma 2.1, there exist $a_1, a_2\in B_{1\over 2}^5$ such that
the maps $h^-_{n,\varepsilon}:V_{n,\epsilon}^-\to
S^4_{1-\epsilon}$, $h^+_{n,\varepsilon}:V_{n,\epsilon}^+\to
S^4_{1+\epsilon}$ defined by
$$h^-_{n,\varepsilon}=(\Pi^-_{a_1}|_{S^4_{1-\varepsilon}})^{-1}\circ
\Pi^-_{a_1}\circ u_n, \quad
h^+_{n,\varepsilon}=(\Pi^+_{a_2}|_{S^4_{1+\varepsilon}})^{-1}\circ
\Pi^+_{a_2}\circ u_n$$ satisfy
$$h^-_{n,\varepsilon} -u_n\in W^{2,2}_0(V^-_{n,\varepsilon},\R^5),
\quad h^+_{n,\varepsilon} -u_n\in W^{2,2}_0
(V^+_{n,\varepsilon},\R^5),$$
$$\int_{V^-_{n,\varepsilon}}|\nabla^2 h^-_{n,\varepsilon}|^2\,dx
\leq C\int_{V^-_{n,\varepsilon}}(|\nabla^2 u_n|^2+|\nabla
u_n|^4)\,dx
$$

and

$$\int_{V^+_{n,\varepsilon}}|\nabla^2 h^+_{n,\varepsilon}|^2\,dx
\leq C\int_{V^+_{n,\varepsilon}}(|\nabla^2 u_n|^2+|\nabla
u_n|^4)\,dx. $$

We now define
$$ w_{n,\varepsilon}(x) =\cases h_n^+(x) &\text { for } x\in
V_{n,\varepsilon}^+\\
h_n^-(x) &\text { for } x\in
V_{n,\varepsilon}^-\\
u_n(x) &\text { for } x\notin V_{n,\varepsilon}^+\cup
V_{n,\varepsilon}^-.
\endcases
$$

Then it is easy to see that $w_{n,\varepsilon}\in
W^{2,2}_\phi(B^5, S^4)$ has only a finitely many singular points
in $\Omega$ and satisfies

$$\int_{V^+_{n,\epsilon}\cup V_{n,\epsilon}^-}|\nabla ^2w_{n,\epsilon} |^2\,dx\leq
C\int_{V^+_{n,\varepsilon}\cup V^-_{n,\varepsilon}}|\nabla^2
u_n|^2+|\nabla u_n|^4)\,dx\to 0, \ \hbox{ as }n\to\infty $$

since $u_n\to u$ strongly in $W^{2,2}$ and
$\lim_{n\rightarrow\infty}|V_{n,\epsilon}^+\cup V_{n,\epsilon}^-|=
0$.

Finally, to obtain the desired approximation, we only have to
project $w_{n,\epsilon}$ on $S^4$ and let $\varepsilon\to 0$.
\qed

\enddemo

\proclaim {Proposition A3} Let $\Phi(x)={x\over |x|}:
B^5\to S^4$. Then there exist infinitely many biharmonic maps
$\{u_i\}\subset W^{2,2}_\Phi(B^5,S^4)$, each of which is smooth
away from a closed set $\Sigma_i$ with ${\Cal
H}^{1-\delta}(\Sigma_i)=0$ for some $\delta>0$.
\endproclaim
\demo {Proof} It is based on some modifications of \cite {BBC}.
First, let $u_0\not=\Phi$ be a given map in $W^{2,2}_{\Phi}(B^5,S^4)$
having finitely many interior singular points.
For $0<\lambda <1$, let $u_{\lambda}\in W^{2,2}_\Phi(B^5,S^4)$
be a minimizer in $W^{2,2}_\Phi(B^5,S^4)$ of
$$\hat{\Bbb H}_{\lambda}(v):=\Bbb H(v)+16\sigma_4\lambda L(u,u_0)$$
where
$$L(u,u_0)=\frac 1{\sigma_4}\sup_{\xi:\Omega\to \R, \|\nabla\xi\|_{L^\infty}\le 1}
\int_\Omega (D(u)\cdot \nabla\xi-D(u_0)\cdot \nabla \xi)\,dx.$$

We remark that theorem B, corollary 3.2, and Lemma 4.3  also hold
for minimizers of $\hat\Bbb H_\lambda$.

We claim: for  $0<\lambda <1$, $u_{\lambda}\not=\Phi$. For,
otherwise, $\Phi$ is a minimizer for both $\Bbb H$ and $\hat{\Bbb
H}_{\lambda}$. In particular, setting $\Phi(t)=\Phi \circ \eta
(t)$, we have
$$\left .\frac d{dt} L(\Phi(t), u_0)\right |_{t=0}=0\tag 5.8$$
where $\eta (t)$ is a smooth family of diffeomorphisms from $B^5$
into itself, satisfying $\eta (0)=Id$ and $\eta (t)=Id$ on
$\partial B$. This is impossible, for we can choose suitable $\eta
(t)$ such that $L(\Phi(t), u_0 )=L(\Phi , u_0)-t+o(t)$ as $t\to
0$.

For a fixed $\lambda_1\in (0,1)$, let
$$A_1=\min \{ \Bbb H(v): \text{ $v$ is a minimizer of $\hat{\Bbb H}_{\lambda_1}$
 }\}.$$
Then there exists a map $u_{\lambda_1}\in W^{2,2}_\Phi(B^5,S^4)$
which minimizes $\hat{\Bbb H}_{\lambda_1}$ such that
$\Bbb H(u_{\lambda_1})=A_1$.
Moreover, since $$\Bbb H (\Phi)<A_1$$
there exists a sufficiently small $0<\lambda_2 <\lambda_1$ such that
$$\Bbb H(\Phi)
+\lambda_2 [\Bbb H (u_0)+\Bbb H(\Phi)] <A_1.
$$
Let $u_{\lambda_2}\in W^{2,2}_\Phi(B^5,S^4)$ be a minimizer of $\hat{\Bbb H}_{\lambda_2}$.
Then we have
$$\split&\Bbb H(u_{\lambda_2})\le
{\hat{\Bbb H}}_{\lambda_2}(u_{\lambda_2})\leq \hat{\Bbb H}_{\lambda_2}(\Phi)\cr
&\leq \Bbb H (\Phi)
+\lambda_2 \int_{B^5}[|\D \Phi|^4+|\D u_0|^4]\,dx\le \Bbb H(\Phi)
+\lambda_2 [\Bbb H(\Phi)+\Bbb H(u_0)]<A_1. \endsplit
$$
This implies that $u_{\lambda_2}$ is different from both $\Phi$ and $u_{\lambda_1}$
Iterating this construction, we find
infinitely many biharmonic maps $u_{\lambda_l}$. By Theorem 3.2
and Lemma 4.3, each $u_{\lambda_l}$ is partially regular. This
proves Proposition A3. \qed
\enddemo

\enddocument

\end